\documentclass[preprint,11pt]{elsarticle}

\usepackage{amssymb}
\usepackage{amsmath}
\usepackage{amsfonts}
\usepackage{amssymb}
\usepackage{indentfirst,latexsym,bm}
\usepackage{amsthm}
\usepackage{color,xcolor}
\usepackage[all]{xy}
\usepackage[colorlinks=false,citecolor=blue]{hyperref}
\usepackage{mathptmx, amsmath, amssymb, amsfonts, amsthm, mathptmx, enumerate, color, mathrsfs}
\input{mathrsfs.sty}
\newtheorem{theorem}{Theorem}[section]
\newtheorem{lemma}[theorem]{Lemma}
\newtheorem{corollary}[theorem]{Corollary}
\newtheorem{proposition}[theorem]{Proposition}

\theoremstyle{definition}

\theoremstyle{definition}
\newtheorem{example}{Example}[section]

\theoremstyle{remark}
\newtheorem{remark}{Remark}[section]
\theoremstyle{question}

\numberwithin{equation}{section}

\journal{XXX}

\begin{document}

\begin{frontmatter}



\title{Some inequalities for adjointable operators on Hilbert $C^*$-modules}
\author[PSUT]{Mohammad Sababheh}
\ead{sababheh@yahoo.com; sababheh@psut.edu.jo}
\author[MB]{Hamid Reza Moradi}
\ead{hrmoradi@mshdiau.ac.ir}
\author[shnu]{Qingxiang Xu}
\ead{qingxiang\_xu@126.com}
\author[shnu]{Shuo Zhao}
\ead{z13427525400@163.com}
\address[PSUT]{Department of Basic Sciences, Princess Sumaya University for Technology, Amman, Jordan}
\address[MB]{Department of Mathematics, Mashhad Branch, Islamic Azad University, Mashhad, Iran}
\address[shnu]{Department of Mathematics, Shanghai Normal University, Shanghai 200234, PR China}
\begin{abstract}The main purpose of this paper is, in the general setting of the adjointable operators on Hilbert $C^*$-modules,  to develop two new tools that can be applied to deal with the positive solutions of certain operator equations, the operator norm as well as the numerical radius, respectively.
Among other things, the positivity of a $2\times 2$ block operator matrix is clarified without any preconditions on its entries, and
a generalized version of the mixed Schwarz inequality with a parameter is  derived. Numerical examples are provided to illustrate the non-triviality of this newly obtained  inequality.
\end{abstract}


\end{frontmatter}



\section{Introduction and the preliminaries}

Given a $C^*$-algebra $\mathfrak{A}$ and  (right) Hilbert $\mathfrak{A}$-modules $H$ and $K$ \cite{Lance}, let
$\mathcal{L}(H, K)$ be the set consisting of all adjointable operators from $H$ to $K$. For each $A\in\mathcal{L}(H, K)$, its adjoint operator, range, and null space are denoted by $A^*$,
$\mathcal{R}(A)$ and $\mathcal{N}(A)$, respectively. Furthermore, the square root of $A^*A$ is denoted by $|A|$. In case $H=K$, $\mathcal{L}(H,K)$ is abbreviated to $\mathcal{L}(H)$, whose unit (namely, the identity operator on $H$) is denoted by $I_H$.  An operator $A\in \mathcal{L}(H)$ is positive, written $A\geq 0$, if $\langle Ax,x\rangle\geq0$ for every $x\in{H }$ \cite[Lemma~4.1]{Lance}. For self-adjoint operators $A, B\in\mathcal{L}(H)$, we say that $B\geq A$ if $B-A \geq 0$. Let $\mathcal{L}(H)_{\mbox{sa}}$  denote the set of all self-adjoint elements in $\mathcal{L}(H)$. In the special case that $H$ is a Hilbert space, we use the notation $\mathbb{B}(H)$ instead of $\mathcal{L}(H)$.

Suppose that $H_1$ and $H_2$ are Hilbert modules over a $C^*$-algebra $\mathfrak{A}$. A Hilbert $\mathfrak{A}$-module can be induced from $H_1$ and $H_2$ as
$$H_1\oplus H_2=\left\{(x,y)^T:x\in H_1, y\in H_2\right\},$$ whose $\mathfrak{A}$-valued inner-product is given  by
$$\left\langle (x_1,y_1)^T, (x_2,y_2)^T\right\rangle=\langle x_1,x_2\rangle+\langle y_1,y_2\rangle$$
for every $x_i\in H_1, y_i\in H_2, i=1,2$. As before, let $\mathcal{L}(H_1\oplus H_2)$ be the set of all adjointable operators on $H_1\oplus H_2$. For each $A\in\mathcal{L}(H_1\oplus H_2)$, it is a matter of routine to verify that $A$ can be represented by $A=(A_{ij})_{1\le i,j\le 2}$,
where
$A_{ij}\in\mathcal{L}(H_j,H_i)$ for $i,j=1,2$.

Let $H$ and $K$ be Hilbert modules over a $C^*$-algebra $\mathfrak{A}$, and  $T$ be an element in $\mathcal{L}(H,K)$.
The Moore-Penrose inverse of $T$ (if it exists), written $T^\dag$, is the
unique element $X\in \mathcal{L}(K,H)$ which satisfies
\begin{equation*} \label{equ:m-p inverse} TXT=T, \quad XTX=X, \quad (TX)^*=TX,\quad (XT)^*=XT.\end{equation*}
It is known that for every  $T\in\mathcal{L}(H,K)$, $T^\dag$ exists if and only if $\mathcal{R}(T)$ is closed in $K$ \cite[Theorem~2.2]{X.S.}.

Unless otherwise specified, throughout the rest of this paper, $\mathbb{C}$ is the complex field, $\mathbb{N}$ is the set of all positive integers, $\mathfrak{A}$ is a $C^*$-algebra, $E$, $H$ and $K$ are non-zero Hilbert $\mathfrak{A}$-modules. Let $T\in\mathcal{L}(H\oplus K)$ be defined by
\begin{equation}\label{bolcked matrix T}T=\left(
                                            \begin{array}{cc}
                                              A & C^* \\
                                              C & B \\
                                            \end{array}
                                          \right),
\end{equation}
where $A\in\mathcal{L}(H)$, $B\in\mathcal{L}(K)$ and $C\in\mathcal{L}(H,K)$. Under the precondition of  the closedness of $\mathcal{R}(A)$, the positivity of $T$ is clarified as follows.

\begin{proposition}\label{prop:xu-sheng} {\rm \cite[Corollary~3.5]{X.S.}} Let $T\in\mathcal{L}(H\oplus K)$ be defined by \eqref{bolcked matrix T} such that $\mathcal{R}(A)$ is closed in $H$.
Then $T\ge 0$ if and only if
$$A\ge 0,\quad C^*=AA^\dag C^*,\quad B-CA^\dag C^*\ge 0.$$
\end{proposition}

A similar characterization for the positivity of $T$ can be provided in the case that $\mathcal{R}(B)$ is closed in $K$. It is notable that in some literatures, Proposition~\ref{prop:xu-sheng} serves as one of the main tools in the study of the positive solutions to certain operator equations (see e.g.\,\cite[Lemma~2.12]{MMX}, \cite[Theorem~2.6]{NCX} and \cite[Theorem~5.5]{XSG}). The main purpose of Section~\ref{sec:positivity}
is to develop a new tool by eliminating any preconditions imposed on $A$ and $B$; see  Theorem~\ref{thm:characterization 2 times 2 positive operator} for the details. As a result, some unified improvements can be made in such field as the  positive solutions to certain operator equations  by replacing Proposition~\ref{prop:xu-sheng} with Theorem~\ref{thm:characterization 2 times 2 positive operator}.

Suppose that $H$ and $K$ are Hilbert spaces and $T\in\mathbb{B}(H,K)\setminus \{0\}$. The  mixed Schwarz inequality is known as
 \begin{equation*}\left |\langle Tx,y\rangle\right |^2\le \langle |T|^{2\alpha} x,x\rangle\cdot \langle |T^*|^{2(1-\alpha)} y,y\rangle \quad (x\in H, y\in K)
\end{equation*}
 for every $\alpha\in [0,1]$,  which can be derived by using the polar decomposition of  $T$. However,
 as shown in  \cite[Example~3.15]{LLX-AIOT}, an adjointable operator acting on a Hilbert $C^*$-module may fail to have the polar decomposition. Fortunately, for every adjointable operator on a Hilbert $C^*$-module, it always has the generalized polar decomposition, which is introduced recently in \cite{ZTX} as follows.

\begin{lemma}\label{lem:ZTX}{\rm \cite[Theorem~3.1]{ZTX}} For every $T\in\mathcal{L}(H,K)$ and $\alpha\in (0, 1)$, there exists a unique $U_\alpha\in\mathcal{L}(H,K)$ such that
\begin{enumerate}
\item[{\rm (i)}] $T=U_\alpha|T|^\alpha$ and $T^*=U_\alpha^*|T^*|^\alpha$;
\item[{\rm (ii)}]$U_\alpha^*U_\alpha=|T|^{2(1-\alpha)}$ and $U_\alpha U_\alpha^*=|T^*|^{2(1-\alpha)}$;
\item[{\rm (iii)}] $U_\alpha|T|^\beta=|T^*|^\beta U_\alpha$\, for every $\beta>0$.
\end{enumerate}
\end{lemma}
It is notable that the above generalized polar decomposition is particularly useful in dealing with some inequalities related to Hilbert $C^*$-modules. Indeed, based on Lemma~\ref{lem:ZTX} a generalized version of the mixed Schwarz inequality with a parameter can be derived in the setting of Hilbert $C^*$-modules; see Theorem~\ref{thm:fg=t} for the details. This gives a new tool in dealing with the operator norm as well as the numerical radius; see examples in Section~\ref{sec:examples}.

The rest of the paper is organized as follows. Section~\ref{sec:positivity} is devoted to the characterizations of the positivity of $2\times 2$ block operator matrices, while Section~\ref{sec:mixed inequality} is concerned with the
generalization of the mixed Schwarz inequality. Two numerical examples are provided in Section~\ref{sec:examples} to show the non-triviality of the parameter involved in the generalized version of the mixed Schwarz inequality.

\section{Characterizations of the positivity of $2\times 2$ block operator
matrices}\label{sec:positivity}

We begin with two useful lemmas.

\begin{lemma}{\rm \cite[Proposition 1.1]{Lance}} If $E$ is a semi-inner product $\mathfrak{A}$-module and $x,y\in E$, then
\begin{equation*}\langle x,y\rangle^* \langle x,y\rangle\le \|\langle x,x\rangle\|\langle y,y\rangle.
\end{equation*}
Hence
\begin{equation}\label{Cauchy-Schwarz inequality}\|\langle x,y\rangle\|\le \|x\|\cdot\|y\|.
\end{equation}
\end{lemma}

\begin{lemma}\label{lem:majorization-01}{\rm\cite[Theorem 2.4]{FMX}} For every $T\in {\mathcal L}(E,K)$ and $S\in {\mathcal L}(H,K)$, the following two statements are equivalent:
\begin{enumerate}
\item[{\rm (i)}] $TT^*\le SS^*$;
\item[{\rm (ii)}] $\Vert T^*x\Vert \le \Vert S^*x\Vert$ for every $x \in K$.
\end{enumerate}
\end{lemma}

Now, we provide the technical result of this section as follows.
\begin{theorem}\label{thm:characterization 2 times 2 positive operator} Let $T\in\mathcal{L}(H\oplus K)$ be defined by \eqref{bolcked matrix T}.
Then the following two statements are equivalent:
\begin{enumerate}
\item[{\rm (i)}] $T\ge 0$;
\item[{\rm (ii)}]$A\ge 0$, $B\ge 0$ and
\begin{equation}\label{positivity for C}\left\|\langle Cu,v\rangle\right\|^2\le \|\langle Au,u\rangle\|\cdot \|\langle Bv,v\rangle\|,\quad \forall\,u\in H, v\in K.
\end{equation}
\end{enumerate}
\end{theorem}
\begin{proof}(i)$\Longrightarrow$(ii). Suppose that $T\ge 0$. Since $T^\frac12\in \mathcal{L}(H\oplus K)$, we have
$$T^\frac12=\left(
              \begin{array}{cc}
                S_{11} & S_{12} \\
                S_{21} & S_{22} \\
              \end{array}
            \right)$$ for
some $S_{11}\in\mathcal{L}(H)$, $S_{12}\in\mathcal{L}(K,H)$, $S_{21}\in\mathcal{L}(H,K)$ and $S_{22}\in\mathcal{L}(K)$.  Let
$$X=(S_{11},S_{12}),\quad Y=(S_{21},S_{22}).$$ Then
$T^\frac12=\left(
             \begin{array}{c}
               X \\
               Y \\
             \end{array}
           \right)$  with  $X\in\mathcal{L}(H\oplus K,H)$ and $Y\in\mathcal{L}(H\oplus K,K)$ such that
$$X^*=\left(
        \begin{array}{c}
          S_{11}^* \\
          S_{12}^* \\
        \end{array}
      \right)\in\mathcal{L}(H,H\oplus K),\quad Y^*=\left(
        \begin{array}{c}
          S_{21}^* \\
          S_{22}^* \\
        \end{array}
      \right)\in\mathcal{L}(K,H\oplus K).$$
Hence,  \begin{align*}T=&\left(
                   \begin{array}{c}
                     X \\
                     Y \\
                   \end{array}
                 \right)\left(
                          \begin{array}{cc}
                            X^* & Y^* \\
                          \end{array}
                        \right)=\left(
                                  \begin{array}{cc}
                                    XX^* & XY^* \\
                                    YX^* & YY^* \\
                                  \end{array}
                                \right).
\end{align*}
It follows that
\begin{equation}\label{expressions of A B C}A=XX^*\ge 0,\quad B=YY^*\ge 0,\quad C=YX^*.\end{equation}
So,   by \eqref{Cauchy-Schwarz inequality} and \eqref{expressions of A B C} we know that for every  $u\in H$ and $v\in K$,
\begin{align*}\left\|\langle Cu,v\rangle\right\|^2=&\|\langle YX^*u,v\rangle\|^2=\|\langle X^*u,Y^*v\rangle\|^2\le \|X^*u\|^2\cdot \|Y^* v\|^2\\
=&\|\langle XX^*u,u\rangle\|\cdot \|\langle YY^*v,v\rangle\|=\|\langle Au,u\rangle\|\cdot \|\langle Bv,v\rangle\|.
\end{align*}
This shows the validity of \eqref{positivity for C}.

(ii)$\Longrightarrow$(i). First, we consider the  case that $A$ is both positive and invertible. Namely, $A$ is positive definite. In this case, we have
\begin{equation}\label{rep of T pd case}T=\left(
                     \begin{array}{cc}
                       I_H & 0 \\
                       CA^{-1} & I_K \\
                     \end{array}
                   \right)\left(
                            \begin{array}{cc}
                              A & 0 \\
                              0 & B-CA^{-1}C^* \\
                            \end{array}
                          \right)\left(
                     \begin{array}{cc}
                       I_H & 0 \\
                       CA^{-1} & I_K \\
                     \end{array}
                   \right)^*.
\end{equation}
Given an arbitrary element $\xi$ in $K$, let
$$u=A^{-1}C^*\xi,\quad v=\xi.$$
Substituting the above $u$ and $v$ into  \eqref{positivity for C} yields
$$\left\|A^{-\frac12}C^*\xi\right\|^4\le \left\|A^{-\frac12}C^*\xi\right\|^2\cdot \|B^\frac12\xi\|^2,$$
which clearly gives
$$ \left\|A^{-\frac12}C^*\xi\right\|\le \|B^\frac12\xi\|.$$
It follows from  Lemma~\ref{lem:majorization-01} that
$$CA^{-1}C^*=(A^{-\frac12}C^*)^*A^{-\frac12}C^*\le (B^\frac12)^*B^\frac12=B.$$
Hence, according to \eqref{rep of T pd case} we arrive at $T\ge 0$.

Next, we consider the case that $A$ is only positive. For each $n\in\mathbb{N}$, let
$$T_n=\left(
        \begin{array}{cc}
          A_n & C^* \\
          C & B \\
        \end{array}
      \right),$$
where $A_n=A+\frac1n I_H$, which is positive definite in $\mathcal{L}(H)$. Furthermore, for every $u\in H$ we have
$\langle A_n u,u\rangle\ge\langle Au,u\rangle$, which leads by \cite[Proposition~1.3.5]{Pedersen} to
$$\left\|\langle A_n u,u\rangle\right\|\ge \|\langle Au,u\rangle\|.$$
So, \eqref{positivity for C} is valid whenever $A$ is replaced with $A_n$ for every $n\in\mathbb{N}$. Hence, $T_n\ge 0$ for all $n\in\mathbb{N}$.
Therefore, $T=\lim\limits_{n\to\infty}T_n\ge 0$.
\end{proof}

\begin{remark}{\rm The special case of the preceding theorem can be found in \cite[Lemma~1]{Kittaneh}, where $H$ and $K$ are Hilbert spaces. See also
\cite{MMX} for other characterizations of the positivity of $T$ in the Hilbert space case.}
\end{remark}

We end this section by providing two propositions concerning the positivity of $T$ defined by \eqref{bolcked matrix T}.
For this, we need a well-known result, whose proof is omitted.
\begin{lemma}\label{lem:norm of a diagonal matrix}
For every $T=\begin{pmatrix}
\begin{array}{cc}
A & 0 \\
0 & B \\
\end{array}\end{pmatrix}\in\mathcal{L}(H\oplus K)$ with $A\in\mathcal{L}(H)$ and $B\in\mathcal{L}(K)$,
we have $\|T\|=\max\big\{\|A\|,\|B\|\big\}$.
\end{lemma}

Our first proposition reads as follows.
\begin{proposition}\label{6} Let $T\in\mathcal{L}(H\oplus K)$ be given by
\eqref{bolcked matrix T} such that $T\ge 0$. Then
\begin{equation}\label{norm ub-half}\|C\|\le \frac12 \|T\|.
\end{equation}
\end{proposition}
\begin{proof}Let $S,U,\widetilde{T}\in\mathcal{L}(H\oplus K)_{\mbox{sa}}$ be defined by
$$S=\left(
      \begin{array}{cc}
        0 & C^* \\
        C & 0 \\
      \end{array}
    \right),\quad U=\left(
                      \begin{array}{cc}
                        I_H & 0 \\
                        0 & -I_K \\
                      \end{array}
                    \right),\quad \widetilde{T}=UTU.
    $$
It is clear that $S^2=\mbox{diag}(C^*C, CC^*)$,  so a direct use of Lemma~\ref{lem:norm of a diagonal matrix} gives
\begin{equation}\label{same norms of S and C}\|S\|=\sqrt{\max\{\|C^*C\|,\|CC^*\|\}}=\|C\|.\end{equation}
Since $U$ is a symmetry (a unitary and self-adjoint element) and $T\ge 0$, we see that $\widetilde{T}\ge 0$ and  $\|\widetilde{T}\|=\|T\|$.
A simple calculation shows that
$$\widetilde{T}=T-2S,\quad T=\widetilde{T}+2S.$$ It follows from the positivity of $T$ and $\widetilde{T}$ that
$$S\le \frac12 T\le\frac12 \|T\|,\quad -S\le\frac12 \widetilde{T}\le\frac12 \|\widetilde{T}\|=\frac12 \|T\|.$$
Hence,
$$-\frac12 \|T\|\le S\le \frac12 \|T\|,$$
which clearly yields $\|S\|\le \frac12 \|T\|$. The conclusion is immediate from \eqref{same norms of S and C}.
\end{proof}

\begin{remark} {\rm For finite-dimensional inner-product spaces $H$ and $K$, more properties can be found in \cite{Tao} when $T$ defined by \eqref{bolcked matrix T} is a positive semi-definite matrix.}
\end{remark}

\begin{proposition}\label{5} Suppose that the operator $T$ given by \eqref{bolcked matrix T} is positive, and there exists a contraction $U\in\mathcal{L}(H,K)$ such that
\begin{equation}C=U|C|,\quad UU^*B=B.\end{equation}
Then
\begin{equation}\|T\|\le \|A+U^*BU\|.\end{equation}
\end{proposition}
\begin{proof}It follows from Theorem~\ref{thm:characterization 2 times 2 positive operator} that both $A$ and $B$ are positive, hence
$$B=B^*=(UU^*B)^*=BUU^*.$$
Since $U$ is a contraction and $C=U|C|$, by \cite[Lemma~5.1]{DXZ} we have  $U^*C=|C|$. Let
\begin{align*}&\widetilde{U}=\left(
                  \begin{array}{cc}
                    I_H & 0 \\
                    0 & U \\
                  \end{array}
                \right),\quad S=\left(
                                  \begin{array}{cc}
                                    A & |C| \\
                                    |C| & U^*BU \\
                                  \end{array}
                                \right),\\
&\widetilde{S}=\left(
                                                                                 \begin{array}{cc}
                                                                                   0 & -I_H \\
                                                                                   I_H & 0 \\
                                                                                 \end{array}
                                                                               \right)S\left(
                                                                                 \begin{array}{cc}
                                                                                   0 & I_H \\
                                                                                   -I_H & 0 \\
                                                                                 \end{array}
                                                                               \right).\end{align*}
Then $\widetilde{U}\in\mathcal{L}(H\oplus H, H\oplus K)$, $S,\widetilde{S}\in\mathcal{L}(H\oplus H)$ and
$$\widetilde{U}^*T\widetilde{U}=S, \quad \widetilde{U}S\widetilde{U}^*=T,\quad \widetilde{S}=\left(
                                                                                               \begin{array}{cc}
                                                                                                 U^*BU & -|C| \\
                                                                                                 -|C| & A \\
                                                                                               \end{array}
                                                                                             \right).
$$
From Lemma~\ref{lem:norm of a diagonal matrix}, we have
$$ \|\widetilde{U}\|^2=\|\widetilde{U}^*\widetilde{U}\|=\max\{1,\|U^*U\|\}=\max\{1,\|U\|^2\}=1,$$
so $\|S\|\le \|\widetilde{U}^*\|\cdot \|T\|\cdot \|\widetilde{U}\|=\|T\|.$
Similar reasoning shows that $\|T\|\le \|S\|$. Consequently, $\|T\|=\|S\|$. By assumption $T\ge 0$, so we have $S\ge 0$, which in turn gives $\widetilde{S}\ge 0$. Therefore,
\begin{equation*}\|T\|=\|S\|\le\|S+\widetilde{S}\|=\|A+U^*BU\|,
\end{equation*}
since $S+\widetilde{S}=\mbox{diag}(A+U^*BU,A+U^*BU)$.
\end{proof}
\begin{remark}{\rm
For more results related to the above proposition in the matrix case, the reader is referred to \cite[Section~2]{moorr}.}
\end{remark}

\section{A generalized version of the mixed Schwarz inequality}\label{sec:mixed inequality}

To avoid triviality, in this section we only consider non-zero operators. For every Hilbert module $V$ over a $C^*$-algebra and $A\in\mathcal{L}(V)\setminus\{0\}$, let $A^0=I_V$. With such a convention, we provide a generalized version of the mixed Schwarz inequality with parameters as follows.

\begin{theorem}\label{thm:fg=t} For every $T\in\mathcal{L}(H,K)\setminus \{0\}$ and $\alpha\in (0,1]$,  we have
\begin{equation}\label{equ:fg=t with inf} \|\langle Tx,y\rangle\|\le \|T\|^{1-\alpha}\cdot \|f(|T|^\alpha) x\|\cdot \|g(|T^*|^\alpha)   y\|,\quad \forall x\in H, y\in K,
\end{equation}
where $f$ and $g$ are non-negative continuous functions on the interval $[0,+\infty)$ such that $f(t)g(t)=t$ for all $t\in [0,+\infty)$.
\end{theorem}
\begin{proof}First, we consider the case that $\alpha\in (0,1)$. In this case, $T=U_\alpha |T|^\alpha$ such that items (i)--(iii) in Lemma~\ref{lem:ZTX} are all satisfied. Since $|T|^\alpha U_\alpha^*=T^*=U_\alpha^*|T^*|^\alpha$, we have
$$h(|T|^\alpha)\cdot U_\alpha^*=U_\alpha^* \cdot h(|T^*|^\alpha)$$
for every  function $h$ that is continuous on $[0,+\infty)$. It follows that for every $x\in H$ and $y\in K$,
\begin{align*}\|\langle Tx,y\rangle\|=&\|\langle U_\alpha |T|^\alpha x,y\rangle\|=\|\langle  |T|^\alpha x, U_\alpha^* y\rangle\|=\|\langle g(|T|^\alpha) f(|T|^\alpha) x, U_\alpha^* y\rangle\|\\
=&\|\langle f(|T|^\alpha) x,  g(|T|^\alpha)\cdot  U_\alpha^* y\rangle\|=\|\langle f(|T|^\alpha) x,  U_\alpha^*\cdot g(|T^*|^\alpha)   y\rangle\|\\
\le& \|f(|T|^\alpha) x\|\cdot \|U_\alpha^*\cdot g(|T^*|^\alpha)   y\|\le \|f(|T|^\alpha) x\|\cdot \|U_\alpha^*\| \cdot \|g(|T^*|^\alpha)   y\|.
\end{align*}
Since \begin{align*}\|U_\alpha^*\|=\sqrt{\|U_\alpha^* U_\alpha\|}=\sqrt{\big\| |T|^{2(1-\alpha)}\big\|}=\big \| |T|\big\|^{1-\alpha}=\|T\|^{1-\alpha},
\end{align*}
this gives
\begin{equation}\label{equ:parameter between 0 and 1}\|\langle Tx,y\rangle\|\le \|T\|^{1-\alpha}\cdot \|f(|T|^\alpha) x\|\cdot \|g(|T^*|^\alpha)   y\|,\quad\forall x\in H,y\in K.\end{equation}

Next, we consider the case that $\alpha=1$. A direct use of \eqref{equ:parameter between 0 and 1} yields
\begin{equation}\label{inq wrt n variable}\|\langle Tx,y\rangle\|\le \|T\|^{1-\alpha_n}\cdot \|f(|T|^{\alpha_n}) x\|\cdot \|g(|T^*|^{\alpha_n})   y\|, \quad\forall x\in H,y\in K,
\end{equation}
where $\alpha_n=\frac{n}{n+1}$ for each $n\in\mathbb{N}$. By assumption $\|T\|>0$, so
 $\lim\limits_{n\to\infty}\|T\|^{1-\alpha_n}=1$.
Meanwhile,
\begin{align*}&\big\| |T|^{\alpha_n}\big\|=\big \| |T|\big \|^{\alpha_n}=\|T\|^{\alpha_n}\le (\|T\|+1)^{\alpha_n}\le \|T\|+1,\quad \forall n\in\mathbb{N}.
\end{align*}
Similar reasoning shows that
\begin{equation*}\big\| |T^*|\big\|^{\alpha_n}\le \|T\|+1,\quad \forall n\in\mathbb{N}.\end{equation*}
Let $M=\|T\|+1$ and
$$h_n(t)=t^{\alpha_n},\quad h(t)=t,\quad \forall\,t\in [0,M].$$
It follows easily from Dini's theorem that
the sequence $\{h_n\}$ converges uniformly to $h$ on $[0,M]$, which means that
$$\lim_{n\to\infty} \big\| |T|^{\alpha_n}-|T|\big\|=0,\quad \lim_{n\to\infty} \big\| |T^*|^{\alpha_n}-|T^*|\big\|=0.$$
Consequently,
$f(|T|^{\alpha_n})\to f(|T|)$ and $g(|T^*|^{\alpha_n})\to g(|T^*|)$ in the norm topology.
Taking the limit as $n\to \infty$, we see from
 \eqref{inq wrt n variable} that
\begin{equation}\label{equ:fg=t} \|\langle Tx,y\rangle\|\le \|f(|T|)x\|\cdot \|g(|T^*|)y\|,\quad \forall x\in H,y\in K.
\end{equation}
In view of \eqref{equ:parameter between 0 and 1} and \eqref{equ:fg=t}, the conclusion follows.
\end{proof}

\begin{remark}{\rm Inequality \eqref{equ:fg=t} is known in the Hilbert space case; see e.g.\,\cite[Theorem~1]{Kittaneh}.}
\end{remark}

As a consequence of the preceding theorem, we have the following mixed Schwarz inequality in the Hilbert $C^*$-module case.

\begin{corollary}\label{cor:generalized mixed inequality} For every $T\in\mathcal{L}(H,K)\setminus \{0\}$ and $\alpha\in [0,1]$,  we have
 \begin{equation}\label{g-mixed inq}\left\|\langle Tx,y\rangle\right\|^2\le \left\|\langle |T|^{2\alpha} x,x\rangle\right\|\cdot \big\|\langle |T^*|^{2(1-\alpha)} y,y\rangle\big\|,
 \quad \forall x\in H, y\in K.
\end{equation}
\end{corollary}
\begin{proof}Given an arbitrary $\alpha\in (0,1)$, let $f(t)=t^{\alpha}$ and $g(t)=t^{1-\alpha}$ for each $t\in [0,+\infty)$. A simple use of \eqref{equ:fg=t} yields the desired inequality.

For each $x\in H$ and $y\in K$, a direct application of \eqref{Cauchy-Schwarz inequality} yields
\begin{align*}\left\|\langle Tx,y\rangle\right\|^2\le \|Tx\|^2\cdot \|y\|^2=\left\|\langle |T|^2 x,x\rangle\right\|\cdot \|\langle y,y\rangle\|,
\end{align*}
This shows the validity of \eqref{g-mixed inq} in the case that $\alpha=1$. Since
\begin{equation}\label{symmetry eqality}\left\|\langle Tx,y\rangle\right\|=\left\|\langle x,T^* y\rangle\right\|=\left\|\langle T^*y,x\rangle\right\|,\end{equation}
we see that \eqref{g-mixed inq} is also true for $\alpha=0$.
\end{proof}

\section{Two numerical examples}\label{sec:examples}
The purpose of this section is to show the meaningfulness of the parameter $\alpha$ involved in \eqref{equ:fg=t with inf}.
For every $T\in\mathcal{L}(H,K)\setminus\{0\}$ and $\alpha\in (0,1]$, by \eqref{equ:fg=t with inf}  we have
\begin{align*}\|T\|=\sup\big\{\|\langle Tx,y\rangle\|: x\in H,y\in K, \|x\|=\|y\|=1\big\}\le \lambda_{\alpha},
\end{align*}
where
$$\lambda_{\alpha}=\|T\|^{1-\alpha}\cdot \|f(|T|^\alpha)\|\cdot \|g(|T^*|^\alpha)\|,$$
in which $f$ and $g$ are as in Theorem~\ref{thm:fg=t}. It is interesting to find out an adjointable $T$, together with functions $f$ and $g$ as in
Theorem~\ref{thm:fg=t} such that $\lambda_{\alpha}<\lambda_1$ for some $\alpha\in (0,1)$. We provide such an example as follows.

\begin{example}{\rm  Let $H=\mathbb{C}^2$, $T=\mbox{diag}(1,1/2)\in\mathcal{L}(H)$, $\alpha=1/2$ and
$$f(t)=\frac{t}{(1-t)^2+1},\quad g(t)=(1-t)^2+1, \quad \forall t\in [0,+\infty).$$
It is clear that $|T|=|T^*|=T$, $\|T\|=1$ and
\begin{align*}&\|f(|T|)\|=\max\big\{f(1),f(1/2)\big\}=1,\\
&\|f(|T|^\alpha)\|=\max\big\{f(1),f\left[(1/2)^\alpha\right]\big\}=1,\\
&\|g(|T^*|)\|=\max\big\{g(1),g(1/2)\big\}=5/4,\\
&\|g(|T^*|^\alpha)\|=\max\big\{g(1),g\left[(1/2)^\alpha\right]\big\}=5/2-\sqrt{2}.
\end{align*}
Hence, $\lambda_{\alpha}=5/2-\sqrt{2}<5/4=\lambda_1$.}
\end{example}

\begin{example}{\rm
Let $H=\mathbb{C}^2, f(t)=t^2+1, g(t)=\frac{t}{t^2+1}, \alpha=\frac{2}{5}, x=\left(
\begin{array}{c}
 1 \\
 2 \\
\end{array}
\right), y=\left(
\begin{array}{c}
 0 \\
 1 \\
\end{array}
\right)$ and $T=\left(
\begin{array}{cc}
 4 & 0 \\
 8 & 3 \\
\end{array}
\right).$
 Numerical calculations show that
\begin{align*}&|\left<Tx,y\right>|=14, \quad \|T\|^{1-\alpha}\approx 3.8228, \quad f(|T|)=T^*T+I_2=\left(
\begin{array}{cc}
 81 & 24 \\
 24 & 10 \\
\end{array}
\right),\\
&g(|T^*|)=|T^*|(TT^*+I_2)^{-1}\approx \left(
\begin{array}{cc}
 0.4213 & -0.1415 \\
 -0.1415 & 0.1692 \\
\end{array}
\right),\\
&|T|^{\alpha}\approx \left(
\begin{array}{cc}
 2.3299 & 0.3752 \\
 0.3752 & 1.2201 \\
\end{array}
\right), \quad |T^*|^{\alpha}\approx \left(
\begin{array}{cc}
 1.3295 & 0.5002 \\
 0.5002 & 2.2205 \\
\end{array}
\right),\\
&f(|T|^{\alpha})\approx \left(
\begin{array}{cc}
 6.5694 & 1.3319 \\
 1.3319 & 2.6293 \\
\end{array}
\right), \quad g(|T^*|^{\alpha})\approx \left(
\begin{array}{cc}
 0.4729 & -0.0549 \\
 -0.0549 & 0.3750\\
\end{array}
\right),\\
&f(|T|)x\approx \left(
\begin{array}{c}
 129 \\
 44\\
\end{array}
\right), \quad g(|T^*|)y\approx \left(
\begin{array}{c}
 -0.1415 \\
 0.1692 \\
\end{array}
\right),\\
&f(|T|^{\alpha})x\approx \left(
\begin{array}{c}
 9.2332 \\
 6.5905\\
\end{array}
\right), \quad g(|T^*|^{\alpha})y\approx \left(
\begin{array}{c}
 -0.0549 \\
 0.3750 \\
\end{array}
\right).
\end{align*}
Consequently, the right-hand side of \eqref{equ:fg=t with inf} becomes
\[ \|T\|^{1-\alpha}\cdot \|f(|T|^\alpha) x\|\cdot \|g(|T^*|^\alpha)   y\|\approx 16.4356,\]
while the right-hand side of \eqref{equ:fg=t} becomes
\[\|f(|T|)x\|\cdot \|g(|T^*|)y\|\approx 30.0631.\]
So, in this example the upper bound in \eqref{equ:fg=t with inf} is sharper than that in \eqref{equ:fg=t}. This also shows the significance of the parameter $\alpha$ involved in \eqref{equ:fg=t with inf}.

 It is notable that there are other examples where \eqref{equ:fg=t} is however better than \eqref{equ:fg=t with inf}. Thus, the upper bound we added in \eqref{equ:fg=t with inf} is a new independent one  that could be useful.}
\end{example}

\end{document}